\documentstyle[12pt,draft]{article}

\input{psfig.sty}

\begin{document}

\title{Interleaver Design for Turbo Codes}
\author{H. R. Sadjadpour\thanks{ \hspace{0.01in} H. Sadjadpour and N. J. A.  Sloane are with
AT\&T Shannon Labs, Florham Park, NJ. E-mail: sadjadpour@att.com and njas@research.att.com.},
N. J. A. Sloane,  
M. Salehi\thanks{ \hspace{0.01in} M. Salehi is with the Department of  Electrical Engineering, Northeastern University .},
 and G. Nebe\thanks{ \hspace{0.01in} G. Nebe is with the Abteilung Reine Mathematik, Universit\"at Ulm, Germany.}}
\date{Nov. 10, 2000}

\maketitle

\begin{abstract} The performance of a Turbo code with short block  length depends critically on the interleaver design. There are two 
major criteria in the design of an interleaver: the distance spectrum of the code and the correlation between
the information input data
and the soft output of 
each decoder corresponding to its parity bits.
This paper describes a new interleaver design  for Turbo codes with short block length based on these two 
criteria. 
A deterministic interleaver suitable for Turbo codes is also described.
Simulation results compare the new interleaver design to  different existing interleavers. 
\end{abstract}


\section{Introduction}
Turbo codes[1] have an impressive near-Shannon-limit error  correcting performance.
The superior performance of Turbo codes 
over convolutional codes is achieved only when the length of  the interleaver is very large, on the order of several thousand bits. 
For large block size interleavers, most random interleavers perform  well.
On the other hand, for some applications it is preferable to have a deterministic interleaver, to reduce the hardware requirements for interleaving and de-interleaving operations. One of the goals of this paper is to propose a  deterministic interleaver design to address this problem.
For short interleavers, the 
performance of the Turbo code with a
random interleaver degrades  substantially up to a point where its bit error rate (BER) performance 
is worse than the BER
performance of convolutional codes with similar computational  complexity. For short block length interleavers, selection of the 
interleaver has a significant effect on the performance of the Turbo  code. In many applications such as voice, delay is an important issue 
in choosing the block size. For these applications, there is a need to  design short block size interleavers that demonstrate acceptable BER 
performance.
Several authors have suggested interleaver designs for  Turbo codes suitable for short block sizes [2-5]. 

There are two major criteria in the design of an interleaver: 1)
the distance spectrum properties (weight distribution) of the code, and 2) 
the correlation between the soft output of each decoder corresponding  to its parity bits and the information input data sequence.
Criterion 2) is sometimes referred to as the iterative decoding  suitability (IDS) criterion [2].
This is a measure of the effectiveness 
of the iterative decoding algorithm and the fact that if these two data  sequences are less correlated, then the performance of the iterative 
decoding algorithm improves. 

The performance of Turbo codes at low BER is mainly dominated by the  minimum effective free distance ($d_{min}$)[13,16]. It has been shown [6] 
that the Turbo code asymptotic performance approaches the $d_{min}$  asymptote. The noise floor that occurs at moderate to high 
signal-to-noise ratios (SNR) is the result of small $d_{min}$ [6]. The  noise floor can be lowered by increasing either the interleaver size or 
$d_{min}$. The latter can be achieved by appropriate choice of  interleaver.
In our approach, maximizing $d_{min}$ is a goal in 
designing the interleaver.

Performance evaluation of Turbo codes is usually based on the  assumption that the receiver is a maximum likelihood (ML) decoder. 
However, Turbo codes actually use a sub-optimal iterative algorithm. A soft output decoding  algorithm such as maximum a posteriori probability (MAP)[7] is used in the iterative algorithm.
The performance of the iterative decoding improves if the information  that is sent to each decoder from the other decoders is less correlated 
with the input information data sequence. Hokfelt et al. [2] proposed the IDS criterion for 
designing an interleaver.
In the interleaver design proposed here, we recommend the use of the IDS criterion with some modifications. 

Trellis termination of Turbo codes is critical, especially when the  interleaver is designed to maximize $d_{min}$.
If this problem is not 
addressed in the design of the interleaver, it can lead to a very small  value for $d_{min}$ because of the existence of data sequences with no trellis 
termination and low output weight, resulting in a degradation in the  performance of the Turbo code.
Papers [8-10] have addressed this question.

The paper is organized as follows.
In section \ref{sect2} random and  S-random interleavers [11] are described. Our approach is based on 
S-random interleavers.
The IDS [2] criterion is also briefly discussed.  In section \ref{sect3}, a two-step S-random interleaver design is presented.
Our approach requires knowing which polynomials are divisible by a primitive
polynomial; this question is addressed in the Appendix.
Section \ref{sect4} describes a deterministic interleaver design based on the results from section \ref{sect3}. 
We conclude the paper by comparing the BER performance of Turbo codes  utilizing our interleaver design 
to other interleavers.

\section{Problem Statements}
\label{sect2}

An interleaver $\pi$ is a permutation 
$i \mapsto \pi(i) $ that changes the order of  a  data sequence
of $N$ input symbols $d_{1},d_{2},\ldots,d_{N}$. If the input data sequence is 
${\bf d}=[d_{1},d_{2},\ldots,d_{N}]$, then
the permuted data sequence is
${\bf d} P$, where $P$ is an interleaving matrix with
a single 1 in each row and
column, all other entries being zero.
Every interleaver  has a corresponding de-interleaver $(\pi^{-1})$ that acts on the 
interleaved data sequence and restores it to its original order.
The de-interleaving matrix is simply the transpose of the interleaving matrix ($P^T$).

A random interleaver is simply a random permutation $\pi$.
For large values of $N$, most random interleavers utilized in Turbo codes perform  well.
However, as the interleaver block size decreases, the 
performance of a Turbo code degrades substantially,
up to a point when its BER performance is worse than that of a convolutional code with 
similar computational complexity.
Thus the design of short interleavers for Turbo codes is an important problem [2-5].

An $S$-random interleaver (where $S=1,2,3, \ldots$) is a
``semi-random'' interleaver constructed as follows.
Each randomly selected
integer is compared with $S$ previously selected
random integers.
If the difference between the current selection and S previous  selections is smaller than S, the random integer is rejected. 
This process is repeated until $N$ distinct integers have been selected.
Computer simulations have shown that if $S 
\leq \sqrt{\frac{N}{2}}$, then this process converges [11] in a  reasonable time. This interleaver design assures that short cycle 
events are avoided.
A short cycle event occurs when two bits are close to  each other both before and after interleaving. 

A new interleaver design was recently proposed based on the performance  of iterative decoding in Turbo codes [2]. Turbo codes utilize an 
iterative decoding process based on the MAP or other algorithms that  can provide a soft output. At each decoding step, some information 
related to the parity bits of one decoder is fed into the other decoder  together with the systematic data sequence and the parity bits 
corresponding to that decoder.
Figure 1 shows this iterative  decoding scheme.
The inputs to each decoder are the input 
data sequence, $d_{k}$, the parity bits $y_{k}^{1}$ or $y_{k}^{2}$, and  the logarithm of the likelihood ratio (LLR) associated with the parity bits from
the other decoder ($W_{k}^{1}$ or $W_{k}^{2}$), which is used as {\em a priori} information. All these inputs are utilized by the decoder to 
create three outputs corresponding to the weighted version of these  inputs. In Figure 1, $\hat{d}_{k}$ represents the weighted version of 
the input data sequence, $d_{k}$. Also $d_{n}$ in the same figure  demonstrates the fact that the input data sequence is fed into the 
second decoder after interleaving. 
The input to each decoder from the other decoder is used as {\em a priori}  information in the next decoding step and corresponds to the 
weighted version of the parity bits. This information will be more  effective in the performance of iterative decoding if it is less 
correlated with the input data sequence (or interleaved input data  sequence). Therefore it is reasonable to
use this as a criterion for 
designing the interleaver. For large block size interleavers, most  random interleavers provide a low correlation between $W_{k}^{i}$ and 
input data sequence, $d_{k}$.
The correlation coefficient,  $r_{W_{k_{1}}^{1},d_{k_{2}}}^{1}$, is defined as the correlation 
between $W_{k_{1}}^{1}$ and $d_{k_{2}}$. It has been shown [2] that  $r_{W_{k_{1}}^{1},d_{k_{2}}}^{1}$ can be analytically approximated by
\begin{equation}
\hat{r}_{W_{k_{1}}^{1},d_{k_{2}}}^{1} = \left\{ \begin{array}{ll} a   \exp^{-c|k_{1} - k_{2}|} & \mbox{if $k_{1} \neq k_{2}$} \\ 0 & \mbox{if 
$k_{1} = k_{2}$} \end{array} \right. 
\end{equation}
where $a$ and $c$ are constants that depend on the encoder feedback  and feedforward polynomials. 
The correlation coefficient at the output 
of the second decoder, 
${\bf r}_{{\bf W}^{2},{\bf d}}^{2}$, is approximated by
\begin{equation}
\hat{\bf r}_{{\bf W}^{2},{\bf d}}^{2} = \frac{1}{2} \hat{\bf r}_{{\bf  W}^{1},{\bf d}}^{1} P (I + \hat{\bf r}_{{\bf W}^{1},{\bf d}}^{1} )
\end{equation}
where the two terms in the right hand side of (2) correspond to the  correlation coefficients between ${\bf W}^{2}$ and the input data, 
i.e., ${\bf d}$  and ${\bf W}^{1}$   [2]. In our notation, $\hat{\bf  r}_{{\bf W}^{2},{\bf d}}^{2}$ represents the correlation coefficient 
matrix and $\hat {r}_{W_{k_{1}}^{2},d_{k_{2}}}^{2}$ represents one  element of this matrix.

Similar correlation coefficients can be computed for the  de-interleaver. The correlation matrix corresponding to de-interleaver, 
$\hat{\bf r'}_{{\bf W}^{2},{\bf d}}^{2}$, is the same as (2) except that
 $P$  is replaced by $P^{T}$.

Then ${\bf V}_{k_{1}}$ is defined to be 
\begin{equation}
{\bf V}_{k_{1}} = \frac{1}{N-1} \sum_{k_{2}=1}^{N} (\hat  {r}_{W_{k_{1}}^{2},d_{k_{2}}}^{2} - 
\bar{\hat {r}}_{W_{k_{1}}^{2},{\bf d}}^{2})^{2} 
\end{equation}
where 
\begin{equation}
\bar{\hat {r}}_{W_{k_{1}}^{2},{\bf d}}^{2} = \frac{1}{N}  \sum_{k_{2}=1}^{N} \hat {r}_{W_{k_{1}}^{2},d_{k_{2}}}^{2} ~.
\end{equation}
$V'_{k_{1}}$ is defined in a similar way using $\hat{\bf r'}_{{\bf  W}^{2},{\bf d}}^{2}$.
The iterative decoding suitability (IDS) measure is then defined as \
\begin{equation}
IDS = \frac{1}{2N} \sum_{k_{1}=1}^{N} ( V_{k_{1}} + V'_{k_{1}})
\end{equation}
A low value of IDS is an indication that the correlation properties between  ${\bf W}^{1}$ and ${\bf d}$ are equally spread along the data sequence 
of length N. An interleaver design based on the IDS condition is  proposed in [12].

\section{Two-step S-random Interleaver Design}
\label{sect3}

A new interleaver design, a two-step S-random interleaver, is presented  here.
The goal is to
increase the minimum effective free distance, $d_{min}$, of the  Turbo code while decreasing or at least not increasing the correlation 
properties between the information input data sequence and $W_{k}^{i}$.
Hokfelt et al. [2,12] introduced the IDS criterion to evaluate the  correlation properties. The two vectors for the computation of IDS in 
(5) are very similar and for most interleavers. Thus it is  sufficient to only use  one of them, i.e., $V_{k_{1}}$. Instead we can 
define a new criterion based on decreasing the correlation coefficients  for the third decoding step, i.e., the correlation coefficients between 
extrinsic information from the second decoder and information input  data sequence. In this regard, the new correlation coefficient matrix, 
$\hat{\bf r'}_{{\bf W}^{2},{\bf d}}^{2}$, is defined as
\begin{eqnarray}
\hat{\bf r'}_{{\bf W}^{2},{\bf d}}^{2} & = &  \frac{1}{2} \hat{\bf  r}_{{\bf W}^{2},{\bf d}}^{2} P^{T} (I + \hat{\bf r}_{{\bf W}^{2},{\bf 
d}}^{2} )  \nonumber \\ 
& = & \frac{1}{4} ( \hat{\bf r}_{{\bf W}^{1},{\bf d}}^{1} + \hat{\bf  r}_{{\bf W}^{1},{\bf d}}^{1}
P \hat{\bf r}_{{\bf W}^{1},{\bf d}}^{1} P^{T} )   \\
& \times & (I + \frac{1}{2} \hat{\bf r}_{{\bf W}^{1},{\bf d}}^{1} P +  \frac{1}{2} \hat{\bf r}_{{\bf W}^{1},{\bf d}}^{1} P \hat{\bf r}_{{\bf 
W}^{1},{\bf d}}^{1} ) \nonumber
\end{eqnarray}
$V_{k_{1}}^{'(new)}$ can now be computed in a similar way to (3) by using  (6). The new iterative decoding suitability ($IDS_{1}$) is then defined 
as 
\begin{equation}
IDS_{1} = \frac{1}{2N} \sum_{k_{1}=1}^{N} (V_{k_{1}} +  V_{k_{1}}^{'(new)})
\end{equation}
A small value for $IDS_{1}$ only guarantees that the correlation  properties are spread equally throughout the data sequence. However, 
this criterion does not attempt to reduce the power of correlation  coefficients, i.e., $(\hat {r}_{W_{k_{1}}^{2},d_{k_{2}}}^{2})^{2}$ and
$(\hat {r'}_{W_{k_{1}}^{2},d_{k_{2}}}^{2})^{2}$.
Therefore, we  recommend
the following additional condition as a second iterative decoding 
suitability criterion:
\begin{equation}
IDS_{2} = \frac{1}{2N^{2}} \sum_{k_{1}=1}^{N} \sum_{k_{2}=1}^{N}  ( (\hat {r}_{W_{k_{1}}^{2},d_{k_{2}}}^{2})^{2} +  (\hat 
{r'}_{W_{k_{1}}^{2},d_{k_{2}}}^{2})^{2} )
\end{equation}
We then use the average of these two values as a new IDS criterion, namely
\begin{equation}
IDS_{(new)} = \frac{1}{2} (IDS_{1} + IDS_{2} )
\end{equation}
Minimizing (9) is then one of our goals in optimizing the interleaver. 

As we described earlier, S-random interleavers avoid short cycle  events. This property guarantees that two bits close to each other 
before interleaving will have a minimum distance of S after  interleaving. More specifically, for information input data $i$ and 
$j$, and permuted data $\pi (i)$ and $\pi (j)$, an S-random  interleaver will guarantee that if $|i-j| \leq S$, then $| \pi (i) - 
\pi (j) | > S$.
However, this does not exclude the possibility that
$ \pi (j) = j$, which
can degrade the  performance of iterative decoding of Turbo codes for this 
particular bit. The larger the distance between $j$ and $\pi (j)$, the  smaller the correlation between the information input data sequence and 
$W_{k}^{i}$. 
We therefore
introduce  an additional measure, $S_{2}$, which is defined to be the minimum 
permissible distance between $j$ and $\pi (j)$ for all  $j=1,2,\ldots,N$. 

Unlike [12], where the interleaver design is based just on the IDS criterion, our  interleaver is designed in two stages. In the first stage, we design an 
interleaver that satisfies the S-random criterion together with the $S_{2}$  condition. In the second stage, we try to increase the minimum 
effective free distance ($d_{min}$) of the Turbo code while considering  the $IDS_{(new)}$ constraint. The design is as follows. We begin by selecting some 
values for $S_1$ and $S_2$.  \\
{\bf Step 1}: Each randomly selected integer $\pi (i)$ is compared with the previous selections 
$\pi (j)$ to check that if
$i - j \leq S_{1}$ then $|\pi (i) - \pi (j)| > S_{1}$.
We also insist that $\pi$ must satisfy
$|i - \pi (i)| > S_{2}$. \\ 
Besides the above conditions, the last $m$ tail bits used for trellis  termination in the first decoder are 
chosen to satisfy $\pi (1) = N$, and if 
$\pi(i) = N - k$ with $k < m$ then $i < N /2$. This condition  will guarantee that trellis termination for the first decoder is 
sufficient and there will not be any low weight sequence at the output  of the second decoder caused by failure of trellis termination.   \\
{\bf Step 2}: Choose the maximum pre-determined weight $w_{det}$ for  input data sequences and the minimum permissible effective free 
distance of the code $d_{min,w_{det}}$. Find all input data sequences of  length N and weight $w_{l} \leq w_{det}$ and their corresponding 
effective free distance $d_{w_{l}}$ for the Turbo encoder with an  interleaver design based on step 1 such that $d_{w_{l}} \leq 
d_{min,w_{det}}$. All these input data sequences are divisible before  and after interleaving by the feedback polynomial (usually a primitive 
polynomial) of the Turbo encoder.
Consider the first input data block of weight $w_{1}$ with non-zero  elements in locations $(i_{1},i_{2},\ldots,i_{w_{1}})$ and 
$d_{min,w_{1}} \leq d_{min,w_{det}}$. Compute $IDS_{(new)}$ based on  (9) for the original interleaver designed in step 1. Set $j=i_{1}+1$ 
and find the pair $(j,\pi (j))$. Interchange the interleaver pairs  $(i_{1},\pi (i_{1}))$ and $(j,\pi (j))$ to create a new interleaver, 
i.e., $(i_{1},\pi (j))$ and 
$(j,\pi (i_{1}))$. Compute the new IDS, $IDS'_{(new)}$, based on the  new interleaver design.
If
$IDS'_{(new)} \leq IDS_{(new)}$,
replace the interleaver by the new one.
Otherwise, set $j=j+1$ and continue.
Repeat this operation for all input data sequences with a  minimum weight of $w_{l} \leq w_{det}$ and $d_{w_{l}} \leq 
d_{min,w_{det}}$.
After completing this operation, return to step 2 and find all input data sequences of weight $w_{l} \leq w_{det}$ 
with $d_{w_{l}} \leq d_{min,w_{det}}$ for the new interleaver. Continue  this step until it converges and there is no input data sequence 
of weight $w_{l} \leq w_{det}$ with $d_{w_{l}} \leq d_{min,w_{det}}$.  Obviously if $d_{min,w_{det}}$ is too large, the second step 
may never converge,
and in this case $d_{min,w_{det}}$  should be reduced.

An interleaver design proposed in [14] and [15] is based on the joint  S-random criteria and elimination of all error patterns of weight 
$w_{i}$.
However, in practice the joint optimization criteria will not converge  easily and therefore the value of $S$ must be reduced and $w_{i}$ 
restricted to only weight two inputs. For weights larger than two, the  convergence of the algorithm is a problem because of the large number of 
possibilities. By separating these two criteria into two steps, we can  easily find the appropriate interleaver satisfying each step 
separately. 
The two steps in the two-step S-random interleaver design are independent operations. The second step tries to increase the minimum effective free distance of the code (based on the interleaver design in the first step) to a pre-determined value ($d_{min,w_{det}}$), while attempting not to increase the correlation between the information input data and the soft output of each decoder corresponding to its parity bits. Obviously, if $d_{min,w_{det}}$ is set to too large a value, the 
second stage of the design may completely change the interleaver produced by the first step and produce an inferior design. This possibility will be illustrated later by simulation.

It is shown in [13] that the feedback polynomials for the recursive systematic convolutional
encoder of Turbo codes should be chosen to be primitive polynomials.
When used for Turbo codes, primitive
polynomials exhibit better spectrum distance properties. The Appendix 
describes how to find all input data sequences of weight  $w_{det}$ that are divisible by a primitive polynomial. This 
information is required for the second step in our approach.

\section{Deterministic Interleaver Design}
\label{sect4}
The following theorem describes a deterministic interleaver based on step 1 in the previous section. 
\subsection*{\bf Theorem 1}
\hspace*{\parindent}
Let $\alpha$ and $N$ be relatively prime natural numbers such that $\alpha -1$ divides $N$,
and let $S_1 = \min \left\{ \alpha, \left\lfloor \frac{N}{\alpha +1} \right\rfloor \right\}$,
$S_2 = \left\lfloor \frac{\alpha -1}{2} \right\rfloor$.
Then there is a permutation $\pi \in S_N$ such that (a)~if $| (i-j) \bmod~N | \le S_1$ and $i \neq j$ then $| (\pi (i) - \pi (j) ) \bmod~N | \ge S_1$, and
(b)~for all $i$, $|(i- \pi (i)) \bmod~N | \ge S_2$.

\noindent{\bf Proof:}
Let $\beta = \lfloor (\alpha -1 ) /2 \rfloor$ and define $\pi: \{1, \ldots, N \} \to \{1, \ldots, N \}$ by $\pi (i) = \alpha i + \beta$, where $\pi (i)$ is to be interpreted as the number
$\pi (i) \in \{1, \ldots, N \}$ that is congruent to $\alpha i + \beta$ modulo $N$.
Since $gcd (\alpha, N ) =1$, $\pi$ is indeed a permutation.
If $\alpha^{-1}$ denotes the inverse of $\alpha~\bmod~N$, then $\pi^{-1} (j) = \alpha^{-1} (j- \beta )$ is the inverse permutation to $\pi$.

(a)~Note that $S_1 \le \alpha$ and $S_1 \le \lfloor N/ (\alpha +1 ) \rfloor$.
Let $i$ and $j$ be elements of $\{1, \ldots, N \}$ with $i \neq j$ and
$| (i-j) \bmod~N | \le S_1$.
Then either
(i)~$1\le i - j \le S_1$ or (ii)~$1 \le N- (i-j) \le S_1$.

In case (i) we have $|(\pi (i) - \pi (j)) \bmod~N | = |\alpha (i-j) \bmod~N | = \min \{\alpha (i-j )$, $N- \alpha (i-j) \}$, and we will show both terms
are $\ge S_1$.
In fact, since $i-j \ge 1$, $\alpha (i-j) \ge \alpha \ge S_1$.
Also, since $i-j \le S_1 \le N/(\alpha +1 )$, we
have $N- \alpha (i-j) \ge N- \alpha N/(\alpha +1) = N/(\alpha +1) \ge S_1$.

In case (ii) we have $1 \le N- (i-j) \le S_1$, so $N-S_1 \le i-j \le N-1$.
But
$$N - \frac{N}{\alpha} \le N - \frac{N}{\alpha +1} \le N- S_1 ~,$$
so $\alpha N - N \le \alpha (i-j) \le \alpha N - \alpha \le \alpha N$,
which means $(\alpha (i-j)$ is trapped between two successive multiples of $N$,
namely
$(\alpha -1 ) N$ and $\alpha N$.
Therefore
\begin{eqnarray*}
|(\pi (i) - \pi (j)) \bmod~N| =
| \alpha (i-j) \bmod~N | 
 = \min
\{ \alpha N - \alpha (i-j), \alpha (i-j) - (\alpha -1)N \}.
\end{eqnarray*}
Again we show both terms are $\ge S_1$.
Since we are in case (ii), $\alpha N - \alpha (i-j) \ge \alpha \ge S_1$.
Secondly, $\alpha (i-j) - (\alpha -1) N \ge \alpha (N-S_1) - (\alpha -1 ) N =
N- \alpha S_1 \ge N- \alpha N / (\alpha +1) = N/ (\alpha +1) \ge S_1$.

(b)~Let $i \in \{1, \ldots, N \}$.
Then $| (i- \pi (i)) \bmod~N | = | (\alpha -1) i + \beta \bmod~N | \,.$
Since $\alpha -1$ divides $N$, and $\beta = \lfloor (\alpha -1 ) /2 \rfloor$,
the last expression is at least $\lfloor (\alpha -1 ) /2 \rfloor = S_2$. \hfill{Q.E.D.} \\

To maximize  the constants $S_1$ and $S_2$, the number $\alpha $
should be close to $\sqrt{N}$.
Then $S_1$ is also about $\sqrt{N}$.
The following elementary consideration shows that one cannot
achieve $S_1 > \sqrt{N}$:
Assume that $S_1 = \sqrt{N}$.
Then the $\sqrt{N} $ values $\pi (1), \ldots , \pi (\sqrt{N})$ 
have pairwise distance $\geq \sqrt{N}$.
Therefore the `balls' with radius $\sqrt{N}/2$  cover the
$\sqrt{N} \sqrt{N} = N $ numbers $\{ 1,\ldots , N\}$ 
completely.
So Theorem 1 yields a solution where $S_1$ is already optimal.  

In some applications such as wireless systems in Rayleigh fading channels, it has been suggested that an additional interleaver be incorporated either before the first encoder or in the path of the systematic data sequence, or alternatively over the entire data sequence (both the systematic data and the parity bits) in order to improve the performance of the system [17]. The deterministic interleaver proposed here can be used for these applications without adding too much complexity to the system.

\section{Simulation Results and Conclusion}
\label{sect5}
 
This section provides simulation results for the BER performance of  Turbo codes
using the new interleaver design and comparisons with
S-random and random interleavers. The constituent encoders are  recursive systematic convolutional codes with 
memory $m=3$ and with
feedback and feedforward generator polynomials $(15)_{oct}$ and  $(17)_{oct}$ respectively. The trellis termination is applied only to 
the first encoder. 

In all the examples, the number of iterations (using the logarithmic version of the BCJR algorithm [7]) is 18. 
For the first two examples, the signal is BPSK with a code rate of  $\frac{1}{3}$. In the first example, the interleaver block size is 192. The BER 
performance of the new interleaver design is compared with S-random and  random interleavers. For the new interleaver,  two interleavers with design parameters  
$(S_{1},S_{2},d_{min,w_{det}},w_{det})=(9,3,20,4)$ and $(9,3,24,4)$ are chosen. For the S-random  interleaver, the value of  $S$ is 9. From Figure \ref{fig:fig2} it can be concluded that the new 
interleaver design performs much better than other interleavers at low BER. It is  also obvious that the error floor for Turbo codes is much lower with the 
new interlearver design because of the larger value of $d_{min}$.
This figure also shows that choosing a very large value for $d_{min,w_{det}}$ can degrade the performance of the Turbo code. For this particular example, the two-step S-random interleaver with $d_{min,w_{det}}=20$ performs better than that with $d_{min,w_{det}}=24$. The appropriate maximum value for $d_{min,w_{det}}$ depends on the length of the interleaver and it is usually obtained by trial and simulations. 
Figure \ref{fig:fig3} compares the BER performance of the two-step S-random interleaver design 
with S-random and random interleavers with a block size of 400. For the 
new interleaver the design parameters are  $(S_{1},S_{2},d_{min,w_{det}},w_{det})=(14,6,26,4)$ and for the 
S-random 
interleaver $S=14$. The two-step S-random interleaver has much
better BER performance than the S-random interleaver at low BER and results in a lower error floor for Turbo codes.
In practice, because the correlation 
properties of the input data and the parity information are decreasing  exponentially, it is  sufficient to choose a small value for $S_2$.

We have also compared the two-step S-random interleaver with Hokfelt's interleaver design.
Hokfelt's approach results in many interleavers for each run of the algorithm with different BER performance. If we choose a random instance of these designs, it may perform worse than the S-random or two-step S-random interleaver design. However, if we choose the best resulting interleaver among them, its performance can be as good as the two-step 
S-random interleaver design. For the interleavers of length 192 and 400 bits, the best interleavers found by Hokfelt's approach can perform as well as the two-step S-random interleavers that were used in examples 1 and 2.

For the last example, the signal is QPSK with a code rate of $\frac{1}{2}$. Equal number of parity bits are punctured from both encoders. The code block length is 1024. Figure \ref{fig:fig4} compares the BER performance of a random interleaver with a deterministic interleaver described in section \ref{sect4} with design parameters $(\alpha, S_{1}, S_{2})=(33,30,16)$, with $\beta$ the same as $S_{2}$. The performance of this deterministic interleaver is slightly worse than that of a random interleaver. However, the interleaving and de-interleaving operations can be carried out algebraically in the receiver and transmitter  thus reducing storage requirements.

\appendix
\section{Polynomials Divisible by a Primitive Polynomial}

Let $R = GF (2) [X]$ be the ring of polynomials with binary  coefficients,
and let $p(X) \in R$ be a primitive irreducible polynomial of degree  $m>1$.
We wish to determine all the polynomials $f(X) \in R$ which have low  weight
and are divisible by $p(X)$.
(The weight of a polynomial is the number of nonzero terms.)

Choose a zero $\alpha$ of $p(X)$.
Then $\alpha $ generates $ GF(2^m )$ as a field.
Since $p(X)$ is primitive, by definition the minimal $n>0$ with
$\alpha ^n =1$ is $n =2^m -1$.
Note that the nonzero elements of $GF(2^m)$ are precisely the
$n$ zeros of the polynomial $X^n-1$.

Since $p(X)$ is irreducible, a polynomial $f(X)  \in R$ is divisible by
$p(X)$ if and only if $f(\alpha ) = 0$.
If $i,j \in \mathbf{N}$ satisfy $i\equiv j $ $(\bmod~n)$, then
$ \alpha ^j = \alpha ^i $, hence $X^i+X^j$ is divisible by $p(X)$.
Let $T_2$ be the set of polynomials $X^i + X^j \in R$ with
$0 \le i < j$, $i\equiv j$ $(\bmod~n)$.
More generally, let $T_{2k}$ $(k=2,3, \ldots )$ be the sum of $k$
disjoint (i.e. all monomials are distinct) terms from $T_2$.

Let $H$ be the Hamming single-error-correcting code with generator  polynomial
$p(X)$, and let $A_w$ be the set of codewords of $H$ of weight $w$,
written in the usual way as polynomials of degree $<n$ corresponding
to residue classes in $R/(X^n -1)$.
Note that $A_i$ is empty unless $i \equiv 3$ or $0 (\bmod~4) $,
i.e., $A_1$, $A_2$, $A_5$, $A_6$, $\ldots $  are empty.

\subsection*{\bf Theorem 2}
\hspace*{\parindent}
Let $f(X) \in R$ have weight $w$ and write
$$f(X)=g(X)+h(X)$$
where $g(X) \in T_{2i}$, $h(X) \in R$ has weight $j$, no two exponents of
$h(X)$ are congruent modulo $n$, and the terms of $g(X)$ and  $h(X)$ are disjoint
(i.e. $w=2i+j$).
Then $f(X)$ is divisible by  $p(X)$ if and only if $\phi (h(X)) \in A_j$
where $\phi$ means ``read exponents $\bmod~n$''.


{\bf Proof:}
\begin{itemize}
\item[``$\Leftarrow$'']
Let $f(X) = g(X) + h(X)$ be as in the theorem. 
Since $\phi(h(X))$ is divisible by $p(X)$, one has $\phi(h)(\alpha ) =
h(\alpha ) = 0$.
Therefore $g(X) \in T_{2i}$ 
and $h(X)$ are both  divisible by $p(X)$ and  so is $f(X)$.
By construction the weight of $f(X)$ is $w=2i+j$.
\item[``$\Rightarrow$'']
Let $f(X) \in R$ be divisible by $p(X)$.
By construction $g(X)$ and hence $h(X)$ is divisible by $p(X)$,
where $\phi(h(X)) \in A_j$ for some $j$. 
Again by construction the weight of $h(X)$ is the weight of $\phi(h(X))$
and the weight of $f(X)$ is $2i+j=w$. \hfill{Q.E.D.}
\end{itemize}

Note that the polynomials $g(X)$ and $h(X)$ are not necessarily unique.
But one may define $g(X)$ by starting from the highest exponent of  $f(X)$ 
and always taking the first term that fits to make the
decomposition unique.

We discuss the first few values of $w$ individually, and illustrate
by taking
$m=3$, $n=7$ and
$p(X) = X^3 + X+1$.
Then $H$ is a Hamming code of length 7, containing 7 words of weight 3,
7 of weight 4, and 1 word of weight 7.

{\bf Weight} $w=1$.
No monomials are divisible by $p(X)$.

{\bf Weight} $w=2$.
A weight two polynomial is divisible by $p(X)$ if and only if it is in  $T_2$.

Examples: $1+X^7$, $X^4 + X^{39}$. \newline
General form: $f(X) = X^i+X^{i+7j}$, $i \geq 0 ,j\geq 1$.

{\bf Weight} $w=3$.
A weight three polynomial is divisible by $p(X)$ if and only if it  reduces
to a weight 3 codeword in $H$ when the exponents are read $\bmod~n$.

Example: The 7 words in $A_3$ are the cyclic shifts of $p(X)$ itself.
So for instance $X^{32} + X^{16} + X^8$ is divisible by $p(X)$,
since it reduces to $X^4 + X^2 + X  = X p (X)\in A_3$. \newline
General form: $f(X) = X^{i+7j}+X^{i+1+7k}+X^{i+3+7l}$, $i,j,k,l \in  \mathbf{Z} $, $i+7j, i+1+7k, i+3+7l \geq 0$.

{\bf Weight} $w=4$.
A polynomial of weight 4 is divisible by $p(X)$ if and only if it is in  $T_4$,
 or it reduces to an element of $A_4$ when the exponents are read  $\bmod~n$.

Examples: $1+X^7 + X^{10} + X^{17} \in T_4$,
$1+ X^2 + X^3 + X^4 \in A_4$.

\bibliographystyle{IEEE}

\begin{thebibliography}{99}
\bibitem{1}
C. Berrou, A. Glavieux, and P. Thitimajshima, {\it "Near  Shannon Limit Error-Correcting Coding and Decoding: Turbo Codes,"} 
Proceeding of IEEE International Conference on Communications '93, Geneva, Volume 2, pp. 1064-1070.
\bibitem{2}
J. Hokfelt, O. Edfors, and T. Maseng, {\it "Turbo Codes:  Correlated Extrinsic Information and its Impact on Iterative Decoding 
Performance,"} Proceeding of IEEE 49th Vehicular Technology Conference '99, Houston, Texas, Volume 3, pp. 1871-1875.
\bibitem{3}
A.K. Khandani, {\it "Group Structure of Turbo Codes with  Applications to the Interleaver Design,"} International Symposium on 
Information Theory,  MIT, Boston, August 1998, pp. 421.
\bibitem{4}
O.Y. Takeshita and D.J. Costello, Jr., {\it "New Classes of  Algebraic Interleavers for Turbo Codes,"} International Symposium on 
Information Theory, MIT, Boston, August 1998, pp. 419.
\bibitem{5}
H. Herzberg, {\it "Multilevel Turbo Coding with Short  Interleavers,"} IEEE Journal on selected areas in Communications, 
vol.16, no. 2, pp. 303-309, February 1998. 
\bibitem{6}
L.C. Perez, J. Seghers, and D.J. Costello, {\it "A Distance  Spectrum Interpretation of Turbo Codes,"} IEEE Trans. on Information 
Theory, vol. 42, No. 6, pp.1698-1709, November 1996.
\bibitem{7}
L. Bahl, J. Cocke, F. Jelinek, and J. Raviv, {\it "Optimum  decoding of linear codes for minimizing symbol error rate,"} IEEE 
Trans. on Inf. Theory, vol. IT-20, pp. 284-287, Mar. 1974.
\bibitem{8}
W. Blackert, E. Hall, and S. Wilson, {\it "Turbo code  termination and interleaver conditions,"} Electronics Letters, Vol. 31, Issue 24,  
pp. 2082-2084, November 23 1995.
\bibitem{9}
A.S. Barbulescu and S.S. Pietrobon, {\it "Terminating the  trellis of Turbo codes in the same state,"} Electronic Letters, vol. 
31, pp. 22-23, January 5 1995.
\bibitem{10}
M.C. Reed and S.S. Pietrobon, {\it "Turbo code termination  schemes and a novel alternative for short frames,"} Seventh IEEE 
International Symposium on Personal, Indoor, and Mobile Communications, Vol. 2,  pp. 354-358, October 15-18 1996.
\bibitem{11}
S Dolinar and D. Divsalar, {\it "Weight Distribution for  Turbo codes Using Random and Nonrandom Permutations,"} JPL Progress 
report 42-122, pp. 56-65, August 15, 1995.
\bibitem{12}
J. Hokfelt, O. Edfors, and T. Maseng, {\it "Interleaver  Design for Turbo Codes Based on the Performance of Iterative 
Decoding,"} Proceeding of IEEE ICC '99, Vancouver, Canada, Vol. 1, pp. 93-97.
\bibitem{13}
S. Benedetto and G. Montorsi, {\it "Design of Parallel  Concatenated Convolutional Codes, "} IEEE Trans. on Comm., vol. 44, no. 
5, pp. 591-600, May 1996.
\bibitem{14}
A.K. Khandani, {\it "Optimization of the interleaver structure for Turbo codes,"} Proceeding of the 1999 Canadian workshop on Information theory, pp. 25-28, June 1999.
\bibitem{15}
J. Yuan, B. Vucetic, and W. Feng, {\it "Combined Turbo Codes  and Interleaver Design,"}  
IEEE Trans. on Comm., vol. 47, no. 4, pp. 484-487, April 1999.
\bibitem{16}
D. Divsalr, and R. J. McEliece, {\it "Effective free distance of Turbo codes,"} Electronics Letters, vol. 32, no. 5, pp. 445-445, February 1996. 
\bibitem{17}
E. k. Hall and S. G. Wilson, {\it "Design and Analysis of Turbo Codes on Rayleigh Fading Channels,"} IEEE Journal on Selected Areas in Communications, vol. 16, no. 2, February 1998. 
\end{thebibliography}

\newpage

\begin{figure}[htb]
\centerline{\psfig{figure=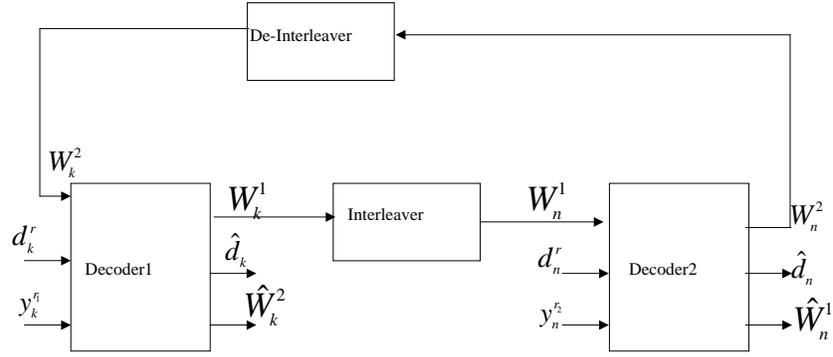,width=5in}}
\caption{Structure of a Tubo decoder.} 
\label{fig:fig1}
\end{figure}

\newpage

\begin{figure}[htb]
\centerline{\psfig{figure=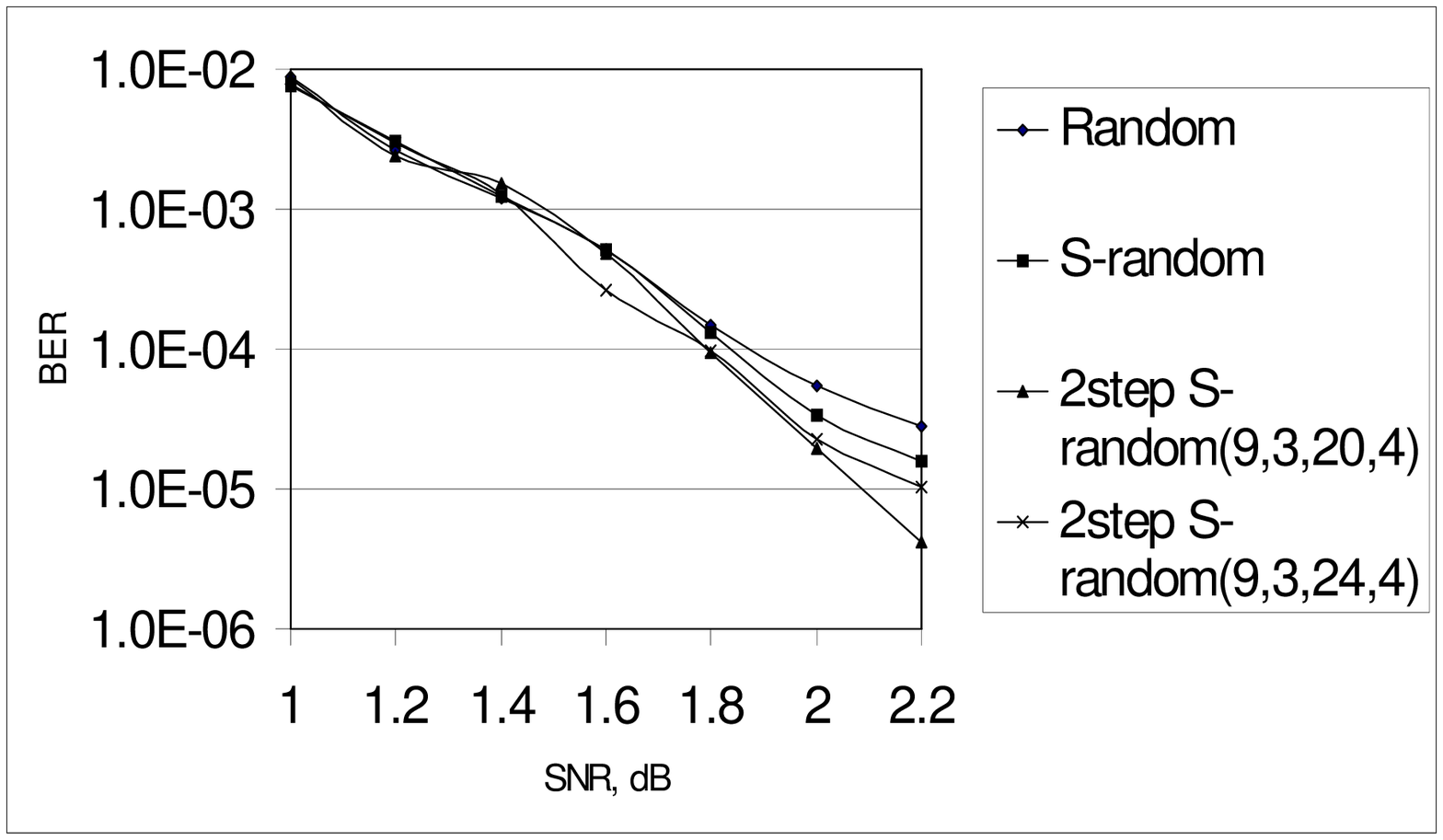,width=5in}}
\caption{Performance of Turbo code for different interleavers of size 192 bits and BPSK signal.} 
\label{fig:fig2}
\end{figure}

\newpage

\begin{figure}[htb]
\centerline{\psfig{figure=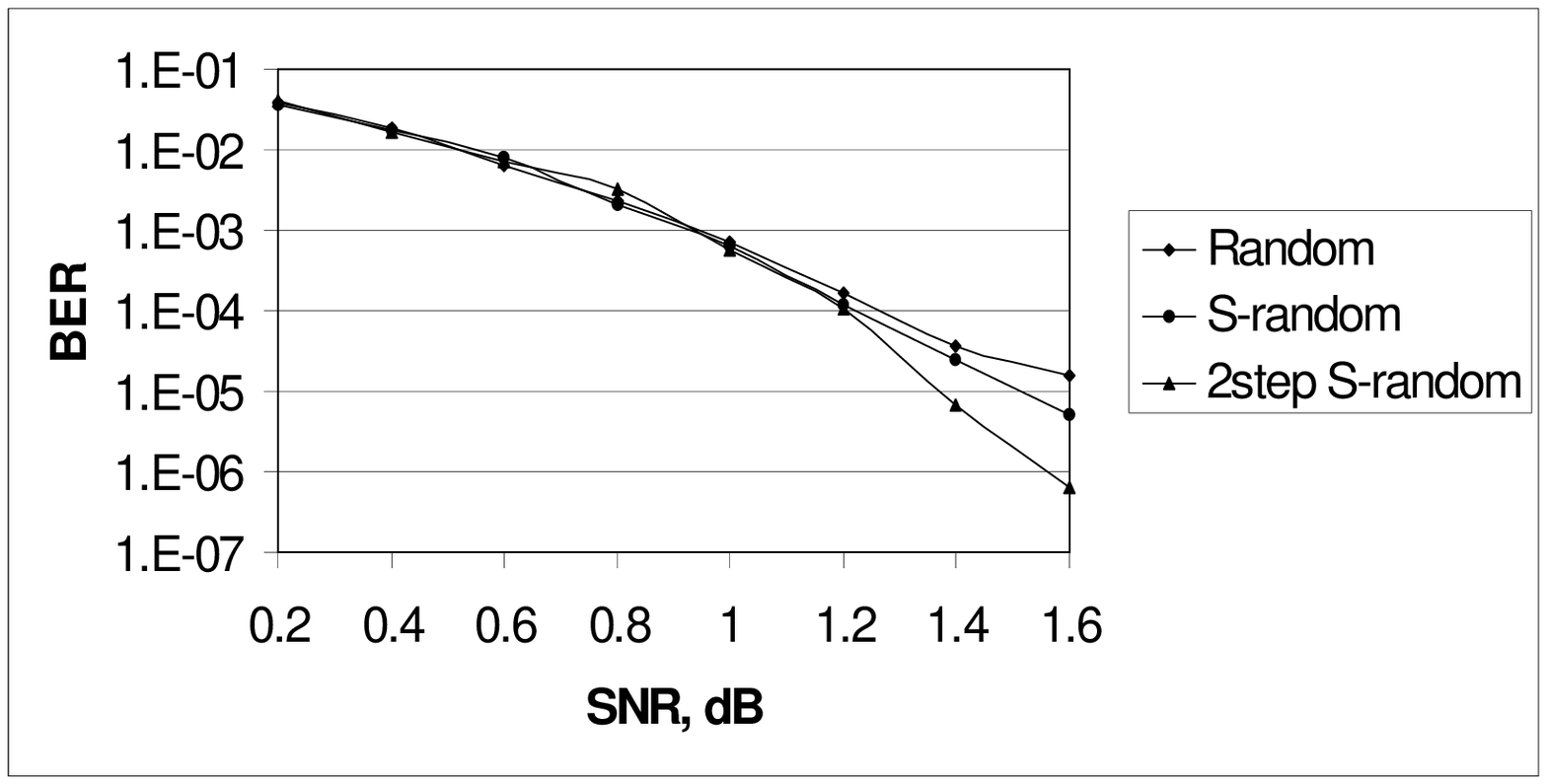,width=5in}}
\caption{Performance of Turbo code for different interleavers of size 400 bits and BPSK signal.} 
\label{fig:fig3}
\end{figure}

\newpage

\begin{figure}[htb]
\centerline{\psfig{figure=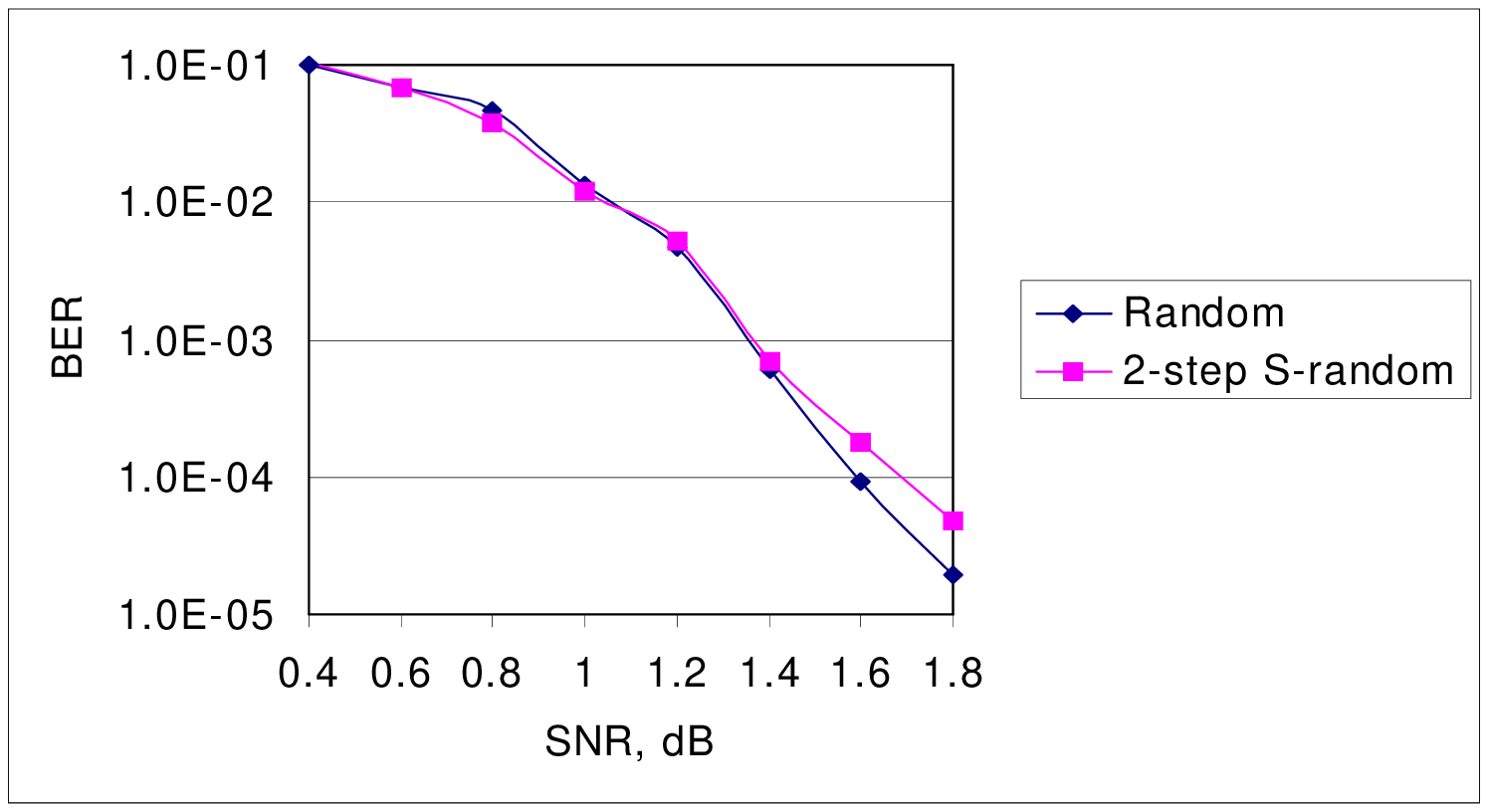,width=5in}}
\caption{Performance of Turbo code for different interleavers of size  1024 bits and QPSK signal.} 
\label{fig:fig4}
\end{figure}

\end{document}